 \documentclass[11pt]{article}
 \usepackage[top=1in,bottom=1in,left=1.5in,right=1.5in]{geometry}
 \usepackage{amsmath}
  \usepackage[dvips]{pstricks} 
  \usepackage{pst-node}
  \usepackage[dvips]{graphicx}
  \usepackage[standard,amsmath,thmmarks]{ntheorem}
  \usepackage{mathrsfs}

  \newcommand{\hs}{\hspace*{\parindent}}
  \newcommand{\per}{\mathop{\mathrm{perm}}\nolimits}
  \newcommand{\permat}{\mathop{\mathrm{perfmat}}\nolimits}

  \newtheorem{theo}{\bfseries \hs Theorem}[section]

  \numberwithin{equation}{section} 

\begin{document}

  \title{The maximum  number of perfect matchings in graphs with a
 given degree sequence}
  \author{
Noga Alon\footnote{School of Mathematics, Tel Aviv
University, Ramat Aviv, Tel Aviv 69978,
Israel, and IAS, Princeton, NJ 08540, USA,
e-mail: nogaa@post.tau.ac.il.
Research supported in part by the Israel Science
Foundation and by a USA-Israeli BSF
grant.} \and
Shmuel Friedland\footnote{Department of Mathematics, Statistics, and Computer Science,
        University of Illinois at Chicago,
        Chicago, Illinois 60607-7045, USA, e-mail friedlan@uic.edu}
  \footnote{
   Visiting Professor, Fall 2007 - Winter 2008,
  Berlin Mathematical School, Berlin, Germany}
  }
  \date{ }

 \maketitle

 \begin{abstract}
We show that the number of perfect matchings in a simple graph
$G$ with an even number of vertices and degree sequence
$d_1,d_2, \ldots ,d_n$ is at most $ \prod_{i=1}^n (d_i
!)^{\frac{1}{2d_i}} $.
 This bound is sharp if and only if $G$ is a union of complete
 balanced bipartite graphs.
     \\[\baselineskip] 2000 Mathematics Subject
     Classification: 05A15, 05C70.
 \end{abstract}


\section{Introduction}

 Let $G=(V,E)$ be an undirected simple graph.
 For a vertex $v\in V$, let $\deg v$ denote  its
 degree.
 Assume that $|V|$ is even, and let $\permat G$
 denote the number of perfect matchings in $G$.  The main result of
 this short note is:
\begin{theo}
\label{maintheo}
\begin{equation}
\label{uppmatbd}
\permat G\le \prod_{v\in V} ((\deg v)!)^{\frac{1}{2\deg v}},
 \end{equation}
 where  $0^{\frac{1}{0}}=0$.  If $G$ has no
 isolated vertices then
 equality holds if and only if $G$ is a disjoint union of
 complete balanced bipartite graphs.
 \end{theo}

 For bipartite graphs the above inequality follows from the
 Bregman-Minc Inequality for  permanents of
$(0,1)$ matrices, mentioned below.

 The inequality (\ref{uppmatbd}) was known to Kahn and Lov\'asz, c.f.
 \cite[(7)]{CK}, but their proof was never published,
and it was recently stated and proved
 independently by the second author in \cite{Fri}.
Here we show that it is a simple consequence of the
Bregman-Minc Inequality.

After our note was published \cite{AlF}, it was pointed out to
us that the inequality that permanent of the adjacency matrix
dominates the square of the number of perfect matching is due
to Gibson \cite{Gib}, and the inequality (\ref{uppmatbd})
appears in \cite[(3.6),p'136]{Ego}.

 \section{The proof}
 Let $A$ be an $n\times n$ $(0,1)$ matrix, i.e.
 $A=[a_{ij}]_{i,j=1}^n \in \{0,1\}^{n\times n}$.  Denote
 $r_i=\sum_{j=1}^n a_{ij}, i=1,\ldots,n$.
 The celebrated Bregman-Minc
 inequality, conjectured by Minc \cite{Min63}
and proved by Bregman \cite{Bre},
 states
 \begin{equation}\label{bminq}
 \per A\le \prod_{i=1}^n (r_i!)^{\frac{1}{r_i}},
 \end{equation}
where equality holds (if no $r_i$ is zero)
iff  up to permutation of rows and columns $A$ is a block diagonal
matrix in which each block is
a square all-$1$ matrix.
\vspace{0.5cm}

\noindent
{\bf Proof of Theorem \ref{maintheo}:}\,
The square of
the number of perfect matchings of $G$ counts ordered pairs
of such matchings. We claim that this is the number of
spanning $2$-regular subgraphs $H $ of $G$
consisting of even cycles (including cycles of length $2$ which are
the same edge taken twice), where each such $H$
is counted $2^s$ times, with $s$ being the number of components (that
is, cycles) of  $H$ with
more than $2$ vertices. Indeed, every union of a pair  of
perfect matchings $M_1,M_2$ is a $2$-regular spanning subgraph $H$
as above,  and for every cycle of length exceeding $2$ in $H$
there are two ways to decide which edges came from
$M_1$ and which from $M_2$.

The permanent of the adjacency matrix $A$ of $G$ also counts
the number of spanning $2$-regular subgraphs $H'$ of $G$, where
now we allow odd cycles and cycles of length $2$ as well. Here,
too, each such $H'$ is counted $2^s$ times, where $s$ is the
number of cycles of $H'$ with more than $2$ vertices, (as there
are $2$ ways to orient each such cycle as a directed cycle and
get a contribution to the permanent).  Thus the square of the
number of perfect matchings is at most the permanent of the
adjacency matrix, and the desired inequality follows from
Bregman-Minc by taking the square root of (\ref{bminq}), where
the numbers $r_i$ are the degrees  of the vertices of $G$.

It is clear that if $G$ is a vertex-disjoint union of balanced complete
bipartite graphs then equality holds in (\ref{uppmatbd}). Conversely,
if $G$ has no isolated vertices and equality holds, then
equality holds in (\ref{bminq}), and no $r_i$ is zero. Therefore,
after permuting the rows and columns of the adjacency matrix of $G$
it is a block diagonal matrix in which every block is an all-1 square matrix,
and as our graph $G$ has no loops, this means that it is
a union of complete balanced bipartite graphs, completing the proof.
$\Box$

\end{document}